%% file: main.tex
\documentclass[runningheads]{llncs}

\usepackage{cite}
\usepackage{amsmath,amssymb,amsfonts}
\usepackage{algorithmic}
\usepackage{graphicx}
\usepackage{textcomp}
\usepackage{xcolor}
\usepackage{multirow}

\usepackage{float}

\usepackage{authblk}

\usepackage[justification=centering]{caption}
\usepackage{subcaption}

\iftrue
\newcommand{\comments}[1]{\footnote{\textcolor{blue}{\textit{#1}}}}

\else
\newcommand{\comments}[1]{}

\fi

\begin{document}

\title{Continuous-time optimal control for trajectory planning under uncertainty}

\author{Ange Valli\orcidID{0000-0001-9483-6834} \and
Shangyuan Zhang\orcidID{0000-0003-0230-8618} \and
Abdel Lisser\orcidID{0000-0003-1318-6679}
}
\authorrunning{A. Valli and al.}
%
\institute{Université Paris-Saclay, CNRS, CentraleSupélec, Laboratoire des signaux et systèmes, 3 Rue Curie Joliot, 91190, Gif-sur-Yvette, France \\
\email{\{ange.valli, shangyuan.zhang, abdel.lisser\}@l2s.centralesupelec.fr}}

\maketitle

\begin{abstract}
This paper presents a continuous-time optimal control framework for generating reference trajectories for autonomous vehicles. The method developed involves chance constraints for modelling uncertainty. A previous work \cite{zhang2022optimal} presented such a model in discrete time and designed for urban driving scenarios only. We obtain better results in continuous time: it models urban driving scenarios with faster computation and better capacity to capture uncertainty as it is less likely to violate the problem's constraints in risky scenarios. It is also robust for optimal control of highway scenarios.
\end{abstract}


\input{Introduction}

\input{Problem_Formulation}

\input{Numerical_Experiments}

\input{Conclusion}

\section*{Acknowledgements}
This research was supported by French government under the France 2030 program, reference ANR-11-IDEX-0003 within the OI H-Code.

\end{document}

%% file: Introduction.tex
\section{Introduction}
\label{section:introduction}

The development of autonomous vehicles is a core research interest in the automotive industry. It aims to meet the challenges of the transport sector. Autonomous vehicle decision-making problems is a vast research field, and trajectory planning is one of its intrinsic components. The trajectory is a spatial-temporal curve connecting the initial position of an ego vehicle (the vehicle on which we perform the control) to a goal to reach, avoiding obstacles and collisions with a target vehicle (another vehicle on which there is no control on its behaviour) \cite{guo2023survey}. The constraints over the dynamics of the vehicle must consider the safety, performance and comfort of the passenger in the proposed solutions. Trajectory planning enhances the safety of autonomous vehicle driving as it predicts complex situations \cite{li2023real} and compensates for human errors of inattention or judgment. The need to model real scenarios and perform trajectory planning on all types of vehicles and robots, including space aircraft \cite{takemura2024uncertainty}, has stimulated research to explore new techniques \cite{le2024stochastic}, including different approaches such as fuzzy-logic systems \cite{hentout2023review} and neural networks \cite{chu2024motion}.

The simulations of autonomous vehicles cannot reproduce the complexity of a real-life environment with all its unexpected events. Furthermore, complex components in autonomous vehicles interact between themselves and each of them can become responsible for a drawback in the nominal functioning. The noise of a sensor, a false measurement due to weather conditions or a technical failure can affect the safeness of the user.

Therefore, a stochastic component in the model is a way for the software to anticipate errors due to hardware components and unpredictable events in the environment. Safety is a significant concern in autonomous vehicles, and the capacity to model mistakes in various driving scenarios addresses this issue. As one cannot design all real-life scenarios for anticipation, simulators must accurately reproduce the most common situations from road traffic. In trajectory planning, a reference trajectory is used as a benchmark to compare with the simulated path of the vehicle.

This paper proposes a continuous-time optimal control problem for generating reference trajectories addressed to autonomous vehicles. The structure of the paper is given as follows: The current Section \ref{section:introduction} presents the background and interest of research. Section \ref{section:related_work} depicts an overview of the literature, tackling the current challenges in trajectory planning and optimal control problems with comparisons of our method with other existing methods, for autonomous vehicles and other types of mobile robots. Section \ref{section:problem_formulation} describes and formulate in equations the continuous-time optimal control problem. We first present the deterministic case and then introduce our stochastic model with chance constraints. Our chance-constrained approach aims at deriving conditions from the probability that constraints are respected. Constraints which are willing to be violated involve a stochastic component. Constraints violations are verified for each time of a solution. Uncertainty is captured by chance constraints derived from the minimal distance between the ego and target vehicles, as the source of stochasticity is the measurement of the position of the target vehicle with respect to the ego vehicle. Section \ref{section:numerical_experiments} is a comparison of different solutions obtained by solving the problem with input samples of scenarios and realisations of the same scenario, including the stochastic component. We assess the models' performances by counting the number of times a constraint is violated in the solution obtained and whether the model can produce a feasible solution or not. It shows the continuous-time model performs better than the discrete-time model, as it computes faster, and the solver finds solutions for trajectory planning over long term horizons, which are helpful in high-speed driving scenarios. Section \ref{section:limitations} discuss limitations based on the study conducted and perspectives of research. Section \ref{section:conclusion} concludes our study and opens it to further research with other approaches than chance constraints and involving other stochastic components to obtain a better model for trajectory planning of autonomous vehicles.

Trajectory planning is a core research area for ensuring safety in autonomous driving. Section \ref{section:related_work} summaries current challenges in optimal control for vehicles and mobile robots. We present our choices to investigate a chance-constrained approach for optimal control of trajectory planning problems and its advantages compared to other methods.

\section{Related Work}
\label{section:related_work}

The trajectory planning problem for autonomous and semi-autonomous vehicles has been studied through various frameworks as the interest in this problem grew within several research fields.

Optimal control with constraints is developed since 2010 \cite{anderson2010optimal} as it enables obtaining a model considering the bounds of the environment and measuring the level of threats depending on the vehicle's current state. Various models were developed for modelling vehicles, among them the unicycle and 4-wheeled models, which are broadly studied nowadays in both robotics \cite{choi2023motion} \cite{sobanski2024predefined} and autonomous vehicles \cite{adzkiya2024control} researches. Extending advances in robotics research to autonomous vehicles opened several aspects of research, such as guidance \cite{degorre2023survey} and trajectory predictions \cite{vishnu2023improving}. Those aspects deal with the challenges of road traffic for the navigation of vehicles in real-life conditions.

In particular, safety constitutes one of the main concerns of intelligent vehicle navigation \cite{zhang2022adversarial}, and new research is carried out. The complexity of the environment of real-life navigation involves considering uncertainties, often treated by including a stochastic component in the modelling. Recent research in the literature \cite{ren2022chance} \cite{vaskov2024friction} \cite{wang2020non} describes it. This line of research is motivated by the desire to have a model that can capture this stochastic component and include it in the vehicle's control to assess the safety level more precisely and address the vehicle's performance. The overview \cite{schwarting2018planning} describes several approaches explored in the literature: game theory, probability, Partially Observable Markov Decision Processes (POMDP) and learning. Various research fields show pros and cons when dealing with autonomous vehicle problems. All techniques present their advantages and weaknesses. Dynamical systems can be modelled by a fractional differential equation, which was already studied for modelling other nonlinear optimal control problems of submarine vehicles \cite{rigatos2020nonlinear} and recently space aircraft \cite{rigatos2024nonlinear}. Those methods are still explored today \cite{zouari2024finite}) and present the advantage of a robust approach with stable Lyapunov control. In those systems, uncertainty can be captured based on neural network approximation \cite{zouari2018neural}. Those methods can be applied to a wide range of systems of different kinds.

Neural networks are used to capture information with the most recent approach based on deep learning \cite{wu2023ccgnet} and deep reinforcement learning \cite{chen2022deep} methods, allowing the neural network to model the uncertainty of the environment and constraints. While this approach captures much information, the computational cost is high. Using deep learning methods to fit real-life contexts is a tough challenge. In autonomous vehicles, embedded systems are used to learn about the environment, but the computation cost can be too high to handle. Therefore, tackling the problem of simulating a driving scenario and performing trajectory planning in real-life conditions should consider those hardware constraints, which can be better handled by other simulation methods. Some methods uses neural networks approximation for motion planning \cite{chu2024motion} or fuzzy-logic systems for capturing uncertainty \cite{hentout2023review}.

Even if those methods present advantages in robustness and stability of system control, our research focuses on optimal control with chance constraints to present a robust model subject to stochastic behaviour. Recent studies were conducted in various fields of research using chance constraints, such as game theory \cite{peng2021chance} \cite{wu2023ccgnet}, linear programming \cite{tassouli2022solving}, geometric optimisation \cite{liu2022distributionally} or decision-making problems modelled by Markov Decision Process (MDP) \cite{varagapriya2023joint}. The advantage of chance constraints is to quantify the level of risk we admit to the decision-maker. Recent research \cite{peng2022bounds} presents advances in joint chance constraint approximations for obtaining bounds to solve problems efficiently based on convex or nonconvex assumptions. In our case, we know there is no guarantee of the solution of the optimal control problem to satisfy all constraints at all time steps because of the stochastic term. Therefore, we capture the probability of satisfaction of a constraint, and we assert this probability to be higher than a threshold $\alpha$ we determine. It provides a guarantee that a certain level of performance can be expected from a model. We show the robustness of the model in our numerical experiments presented in Section \ref{section:numerical_experiments}. These experiments rely on professional software used in the automotive industry, and our method shows its effectiveness in this use case.

Optimisation methods with chance constraints are advantageous as the model is sufficiently robust for giving a feasible solution and respecting a good level of approximation. Recent research \cite{ren2022chance} \cite{zhang2022optimal} presents advances in modelling uncertainty using this framework. The article \cite{ren2022chance} proposes a Gaussian mixture model for uncertainty, allowing convexity properties over the chance constraints to ensure tractability. The authors present simulations with multimodal uncertain obstacles for the trajectory planning of the vehicle, with general polyhedral geometric forms. Previous work in \cite{zhang2022optimal} deals with a constrained nonlinear optimisation problem for generating reference trajectories, with scenarios including an ego vehicle to control and a third-party target vehicle. In both works, the optimisation problems are formulated in discrete time.

In most cases, as stated in the introduction of \cite{ornelas2010discrete}, discrete-time models are derived from approximations of continuous-time models. The interest in tackling discrete-time optimal control problems is to adapt to real-world applications, where numerical sensors are widely used compared to analog sensors. Furthermore, some stability results are guaranteed in discrete time and do not hold for continuous-time edge cases, so it guarantees better stability in the control of systems. Ensuring the stability of some continuous-time frameworks is challenging. Still, they are more representative of real-life settings than discrete-time models, so the choice of model should be addressed in relation to the context of the problem being tackled.

Scientific literature on optimal control problems presents several approaches in the recent research: fractional-order systems, fuzzy-logic systems, neural network approximations, Gaussian mixture models and optimal control problems formulated as nonlinear optimisation problems in discrete time.\\
We focused our research on challenging a Model Predictive Control approach in continuous time by introducing chance constraints to deal with stochastic components coming from measurement errors. Compared to the aforementioned methods, our method presents the advantage of requiring only classical mathematical tools and less computational power than the one required for implementing neural networks or noninteger derivatives. Section \ref{section:problem_formulation} presents the mathematical formulation for our problem and conditions for dealing with stochasticity

%% file: Problem_Formulation.tex
\section{Problem Formulation}
\label{section:problem_formulation}

We divided this section into multiple subsections to define the problem entirely. Subsection \ref{subsection:driving_scenario} defines what a scenario is, with details on inputs that intervene in formulating the optimal control problem. Then, the general form of the optimal control problem is given in Subsection \ref{subsection:general_optimal_control}, and the trajectory planning problem is formulated in Subsection \ref{subsection:reference_trajectory_generator}. Finally, we present the stochastic model in Subsection \ref{subsection:stochastic_model} by introducing chance constraints on the minimum distance between ego and target vehicles as modelling errors of measurements of the position of the target vehicle.

\subsection{Driving scenario}
\label{subsection:driving_scenario}

We present a generalisation in continuous time from the previous work in \cite{zhang2022optimal}. We consider driving scenarios providing information about the trajectory of the vehicles and their environments. This includes other vehicles on the roads and potential obstacles, such as pedestrians. All driving scenarios are defined for a given time horizon.

The prerequisites we consider to be known from the environment to define a driving scenario are the following :

\begin{itemize}
    \item As in discrete time, if we consider $n$ vehicles on the road, for $n \in \mathbb{N}$, $i \in [1, n]$, the trajectory of the $i^{\text{th}}$ vehicle at time $t \in \mathbb{R}$ is represented by its Cartesian coordinates $X_i(t), Y_i(t)$.
    \item The centre lane of the road is a continuous curve represented by its Cartesian coordinates $(x, C(x))_{x \in \mathbb{R^+}}$.
    \item Likewise, the boundaries of the road are represented by $(x, B(x))_{x \in \mathbb{R^+}}$.
    \item We define a maximum speed allowed by regulations on the road, denoted by $v_{max}$.
\end{itemize}

Those elements are the minimum amount of information needed to define a driving scenario. Nonetheless, the reference trajectory needs slightly more input. As in the discrete-time case in \cite{zhang2022optimal}, some additional constraints are also considered:

\begin{itemize}
    \item The initial state of the ego vehicle to control: $z_0$.
    \item A predefined waypoint for anticipating the next actions to undertake. This can take into account lane change, overtaking or steady driving. Our study considers waypoints defined from the centre lane only for simplicity. See Remark \ref{remark:obj_1}.
    \item An optimality criterion to define, here represented by a cost function.   
    \item Additional constraints of the vehicle for considering the passenger's comfort. Those constraints are both cinematic and dynamic.
\end{itemize}

With this scenario definition, we can now formulate the general form of optimal control problem in continuous time with constraints before modelling the reference trajectory generator.

\subsection{Optimal control problem}
\label{subsection:general_optimal_control}
    
Let's recall the formulation of the optimal control problem as in \cite{zhang2022optimal} :

\begin{itemize}
    \item $z(t)$ the control states variable, with $z(t) \in \mathcal{Z}$ and $\mathcal{Z}$ the feasible set of states. $z_{\textrm{init}}$ and $z_{\textrm{term}}$ are initial and terminal states of the system.
    \item $u(t)$ the optimal control input, with $u(t) \in \mathcal{U}$ and $\mathcal{U}$ the feasible set of control inputs.
    \item $t_0$ and $t_n$ are initial and final time, $t_0, t_n \geq 0$.
    \item  $\ell(\cdot)$ the objective function to minimise, $\ell: \mathcal{Z} \times \mathcal{U} \longmapsto \mathbb{R}^+$.
    \item $f(\cdot)$ the function designing the system dynamics of the control state $z(t)$,\\$f: \mathcal{Z} \times \mathcal{U} \longmapsto \mathcal{Z}$.
    \item $c(\cdot)$ is the inequality constraint function, $c: \mathcal{Z} \times \mathcal{U} \longmapsto \mathbb{R}$.
\end{itemize}

The optimal control problem is then given by :

\begin{align}
\underset{z(\cdot), u(\cdot)}{\textrm{min}} \quad & \int_{t_{0}}^{t_{n}} \ell(z(t),\ u(t)) dt \label{eq:general_continuous}\\
\textrm{s.t.} \quad &  \dot{z}(t)=f(z(t), u(t)), \tag{\ref{eq:general_continuous}a} \label{eq:general_control_state}\\
& c(z(t), u(t)) \leq 0, \tag{\ref{eq:general_continuous}b} \label{eq:general_constraints}\\
& z\left(t_{0}\right)=z_{\textrm{init}},\quad z\left(t_{n}\right)=z_{\textrm{term}}, \tag{\ref{eq:general_continuous}c} \label{eq:general_bounds}\\
& z(t) \in \mathcal{Z}, \quad u(t) \in \mathcal{U} \nonumber
\end{align}

With the following constraints on the optimal control problem:\\
\begin{itemize}
    \item The constraint (\ref{eq:general_control_state}) is the control-state equation of the system. It shows the relationship between the state of the system $z(t)$, the control $u(t)$ and the future state of the system.
    \item The constraint (\ref{eq:general_constraints}) describes the boundaries of the possible values over the state and the control of the system. It preserves the stability of our system and considers constraints from the physical problem.
    \item The constraint (\ref{eq:general_bounds}) describes the initial and terminal states of the system. The solution to the optimal control problem needs to give a path to connect them.
\end{itemize}

In this study, our goal is to show the performances of dynamic simultaneous control on the reference trajectory generation, compared to the discretised version of the problem, which is solved as a constrained nonlinear optimisation problem.

We now have all the ingredients to define the optimal control problem of trajectory planning. Subsection \ref{subsection:reference_trajectory_generator} presents the mathematical formulation and the study framework of our research.

\subsection{Modelling a reference trajectory generator}
\label{subsection:reference_trajectory_generator}

The choices for modelling the vehicle align with the discrete-time problem tackled in \cite{zhang2022optimal}. The unicycle kinematic model gives the following state of the ego vehicle at time $t$:
$$z_t = [x_t, y_t, \theta_{t}, v_{t}]^T$$
where $x_t$ is the longitudinal position, $y_t$ is the lateral position, $\theta_t$ is the heading angle of the vehicle and $v_t$ is the linear speed. The control input at time $t$ is given by
$$u_t = [a_t, \omega_t]$$
where $a_t$ is the linear acceleration and $\omega_t$ is the angular velocity. The ego vehicle's control-state relationship is given by :
\begin{gather}
\begin{aligned}
\label{eq:control_state}
    \frac{dz_t}{dt} = f(z_t, u_t)
\end{aligned}
\end{gather}

where $f(z_t, u_t) = [v_t \cos(\theta_t), v_t \sin(\theta_t), \omega_t, a_t]^T$. The optimal control problem proposed for the reference trajectory generation is the following :

\begin{align}
    \min_{\mathbf{u}, \mathbf{z}} \quad & \int_{0}^{T} \mathbf{w}_g * D^2_t(x_t, y_t) + \mathbf{w}_v * (v_r - v_t)^2 + \mathbf{w}_a * a^2_t \nonumber \label{eq:obj_1}\\
    & + \mathbf{w}_{\omega} * \omega^2_t + \mathbf{w}_j * \biggl( \frac{da_t}{dt} \biggl)^2 + \mathbf{w}_h * H(\theta_t)^2 \nonumber \\
    & + \mathbf{w}_p * P(x^{tgt}_t, y^{tgt}_t, x_t, y_t) \hspace{.1cm} dt \\
    \textrm{s.t.} \quad
& \frac{dz_t}{dt} = f(z_t, u_t), \tag{\ref{eq:obj_1}a}\label{eq:ct1}\\
& L(x_t, y_t) \quad \leq \quad 0 , \tag{\ref{eq:obj_1}b}\label{eq:ct2}\\ 
& |v_t|  \leq v_{max} , \tag{\ref{eq:obj_1}c}\label{eq:ct3}\\
& |\omega_t|  \leq \omega_{max} , \tag{\ref{eq:obj_1}d}\label{eq:ct4}\\
& |a_t|  \leq a_{max} , \tag{\ref{eq:obj_1}e}\label{eq:ct5}\\
& \biggl| \frac{da_t}{dt} \biggl| \leq j_{max}, \tag{\ref{eq:obj_1}f}\label{eq:ct6}\\
& K(x^{tgt}_t, y^{tgt}_t, x_t, y_t) \geq d_{min} \tag{\ref{eq:obj_1}g}\label{eq:ct7}\\
& x^{tgt}_t, y^{tgt}_t, x_t, y_t, v_t \in \mathbb{R}^+  \nonumber\\
& a_t \in \mathbb{R} \quad \theta_t, \omega_t \in [-\pi, \pi] \nonumber
\end{align}

With $\mathbf{u}$ and $\mathbf{z}$ vectors representing the input control and the system's state, respectively.

The objective function (\ref{eq:obj_1}) comprises a sum of weighted quantities with terms to control the ego vehicle. The weights are chosen as a trade-off between the comfort and security of the passenger and speedness of the vehicle. Those weights represent the design parameters of the model, and their resolution is discussed in Section \ref{subsection:design_weights}.

\begin{itemize}
    \item $D_t^2(x_t, y_t)$ is the distance to the next waypoint at instant $t$. This term is responsible for controlling the vehicle on the centre line. We have $D_t: \mathbb{R}^2_+ \longmapsto \mathbb{R}^+$.
    \item $(v_r - v_t)^2$ is the $L^2$-distance to the recommended linear speed value $v_r$, considered as a known input parameter.
    \item Minimising quantities $a^2_t$ and $\omega^2_t$ is for penalising input control on linear and angular accelerations, so the values do not become too large, which can make the control unstable.
    \label{jerk_first_mention}
    \item $\bigl( \frac{da_t}{dt} \bigl)^2$ is the jerk term. The jerk is the linear acceleration rate of change over time. A high value of jerk can have physiological effects on the human body. Minimising this quantity allows the vehicle to perform smooth accelerations, which does not make the journey uncomfortable for the passenger. We assume the control variable $a_t$ to be differentiable on $[0, T]$.
    \item $H(\theta_t)^2$ is the distance between the heading angle of the vehicle and the degree of curvature of the centre lane. We have $H: [-\pi, \pi] \longmapsto \mathbb{R}^+$.
    \item $P(x^{tgt}_t, y^{tgt}_t, x_t, y_t)$ is a potential field function useful to regulate the speed of the ego vehicle towards its distance to the target vehicle. This term models the vehicle's Adaptive Cruise Control (ACC) \cite{liu2017path} feature so the ego vehicle can brake quickly if the target vehicle in front is too close. We have $P: \mathbb{R}^4_+ \longmapsto \mathbb{R}^+$.
\end{itemize}

\begin{remark}
\label{remark:obj_1}
    The waypoints $(x^{waypoint}_t, y^{waypoint}_t, \theta^{waypoint}_t)_{t \geq 0}$ are derived from the center lane coordinates directly. Waypoint planning challenges are not tackled in our study.
\end{remark}

The constraints cover several aspects of the modelling:
\begin{itemize}
    \item The constraint (\ref{eq:ct1}) is the control-state relationship (\ref{eq:control_state}) derived from the unicycle kinematic model.
    \item $L(x_t, y_t)$ is the distance between the coordinates of the ego vehicle and the limits of the road defined in constraint (\ref{eq:ct2}). It guarantees the vehicle will stay within the boundaries of the environment. We have $L: \mathbb{R}^2_+ \longmapsto \mathbb{R}^+$.
    \item The constraint (\ref{eq:ct3}) is the maximum linear speed limit of the vehicle.
    \item The constraint (\ref{eq:ct4}) is the maximum angular speed limit of the vehicle.
    \item The constraint (\ref{eq:ct5}) is the maximum linear acceleration limit of the vehicle.
    \item The constraint (\ref{eq:ct6}) is the maximum jerk limit of the vehicle. A large value can be harmful for both the passengers and the vehicle.
    \item $K(x^{tgt}_t, y^{tgt}_t, x_t, y_t)$ is the distance between ego and target vehicles defined in constraint (\ref{eq:ct7}). It guarantees a minimum distance $d_{min}$ between vehicles to prevent collisions. We have $K: \mathbb{R}^4_+ \longmapsto \mathbb{R}^+$.
\end{itemize}

This deterministic formulation of the optimal control problem does not consider the randomness in target vehicle positions $(x^{tgt}_t,y^{tgt}_t)_{t\in\mathbb{R}^+}$. We presented the first contribution of the article, which is the continuous-time formulation of the problem, and Subsection \ref{subsection:stochastic_model} presents the extension to a stochastic model.
\subsection{Stochastic model}
\label{subsection:stochastic_model}
Our stochastic model corresponds to its discrete-time equivalent \cite{zhang2022optimal}. Let's suppose $(x^{tgt}_t)_{t \in \mathbb{R}^+}$ and $(y^{tgt}_t)_{t \in \mathbb{R}^+}$ are sampled from Gaussian processes such as for all $t\in \mathbb{R}^+$:
\begin{align}
    & x^{tgt}_t \sim \mathcal{N}(\mu_{x_t},\sigma_{x_t})
    \label{eq:x_tgt_random}\\ \nonumber \\
    & y^{tgt}_t \sim \mathcal{N}(\mu_{y_t}, \sigma_{y_t})
    \label{eq:y_tgt_random}
\end{align}

Then, the results from \cite{zhang2022optimal} holds. Constraint (\ref{eq:proba_alpha}) is the stochastic constraint, we derive deterministic equivalent second-order conic constraints (\ref{eq:ct8}) and (\ref{eq:ct9}). :
\begin{align}
    & \forall t \in \mathbb{R}^+, \ \mathbb{P}(|x^{tgt}_t - x_t + y^{tgt}_t - y_t| \geq d_{min} ) \geq \alpha \implies
    \label{eq:proba_alpha}\\
    & \forall t \in \mathbb{R}^+, \ x_t + y_t \leq \mu_{x_t} + \mu_{y_t} + d_{min} + \sqrt{\sigma_{x_t}^2 + \sigma_{y_t}^2} \cdot F_N^{-1} (\alpha/2) \label{eq:ct8}\\
    & \forall t \in \mathbb{R}^+, \ x_t + y_t \geq \mu_{x_t} + \mu_{y_t} - d_{min} + \sqrt{\sigma_{x_t}^2 + \sigma_{y_t}^2} \cdot F_N^{-1} (1-\alpha/2) \label{eq:ct9}
\end{align}

We derive our stochastic model from the deterministic model with a modified constraint (\ref{eq:ct7}). This contribution is the main result of this paper. The measurement errors are taken into account in the mathematical formulation, and we include those constraints to guarantee the condition (\ref{eq:proba_alpha}) holds for all scenarios and for any time step $t$. Thanks to Prekopa's work \cite{prekopa2013stochastic}, we can obtain a deterministic equivalent formulation of the constraint. This result is more robust than other methods consisting of the relaxation and approximation of constraints, which give only bounds of the probability, which is not an exact result. In our case, we obtain the equivalence. \\Let $K(x_1,x_2,x_3,x_4) = |x_1 - x_2 + x_3 - x_4|$. Therefore, the constraint (\ref{eq:ct7}) is replaced by the deterministic equivalent constraints (\ref{eq:ct8}) and (\ref{eq:ct9}). The main contribution of this paper is the continuous-time model with constraints (\ref{eq:ct8}) and (\ref{eq:ct9}) defined by the following equations :

\begin{align}
\label{eq:stochastic_model}
    \min_{\mathbf{u}, \mathbf{z}} \quad & \int_{0}^{T} \mathbf{w}_g * D^2_t(x_t, y_t) + \mathbf{w}_v * (v_r - v_t)^2 + \mathbf{w}_a * a^2_t \nonumber \\
    & + \mathbf{w}_{\omega} * \omega^2_t + \mathbf{w}_j * \biggl( \frac{da_t}{dt} \biggl)^2 + \mathbf{w}_h * H(\theta_t)^2 \nonumber \\
    & + \mathbf{w}_p * P(x^{tgt}_t, y^{tgt}_t, x_t, y_t) \hspace{.1cm} dt \\
    \textrm{s.t.} \quad
& \frac{dz_t}{dt} = f(z_t, u_t), \tag{\ref{eq:stochastic_model}a}\\
& L(x_t, y_t) \quad \leq \quad 0 , \tag{\ref{eq:stochastic_model}b}\\ 
& |v_t|  \leq v_{max} , \tag{\ref{eq:stochastic_model}c}\\
& |\omega_t|  \leq \omega_{max} , \tag{\ref{eq:stochastic_model}d}\\
& |a_t|  \leq a_{max} , \tag{\ref{eq:stochastic_model}e}\\
& \biggl| \frac{da_t}{dt} \biggl| \leq j_{max},\tag{\ref{eq:stochastic_model}f}\\
& \forall_{t\in \mathbb{R}^+}, \ x_t + y_t \leq \mu_{x_t} + \mu_{y_t} + d_{min} + \sqrt{\sigma_{x_t}^2 + \sigma_{y_t}^2} \cdot F_N^{-1} (\alpha/2) \tag{\ref{eq:stochastic_model}g}\\
& \forall_{t\in \mathbb{R}^+}, \ x_t + y_t \geq \mu_{x_t} + \mu_{y_t} - d_{min} + \sqrt{\sigma_{x_t}^2 + \sigma_{y_t}^2} \cdot F_N^{-1} (1-\alpha/2) \tag{\ref{eq:stochastic_model}h}\\
& x^{tgt}_t, y^{tgt}_t, x_t, y_t, v_t \in \mathbb{R}^+  \nonumber\\
& a_t \in \mathbb{R} \quad \theta_t, \omega_t \in [-\pi, \pi] \nonumber
\end{align}

\begin{figure}[htp]
    \color{red}
    \centering
    \scalebox{.35}{\includegraphics{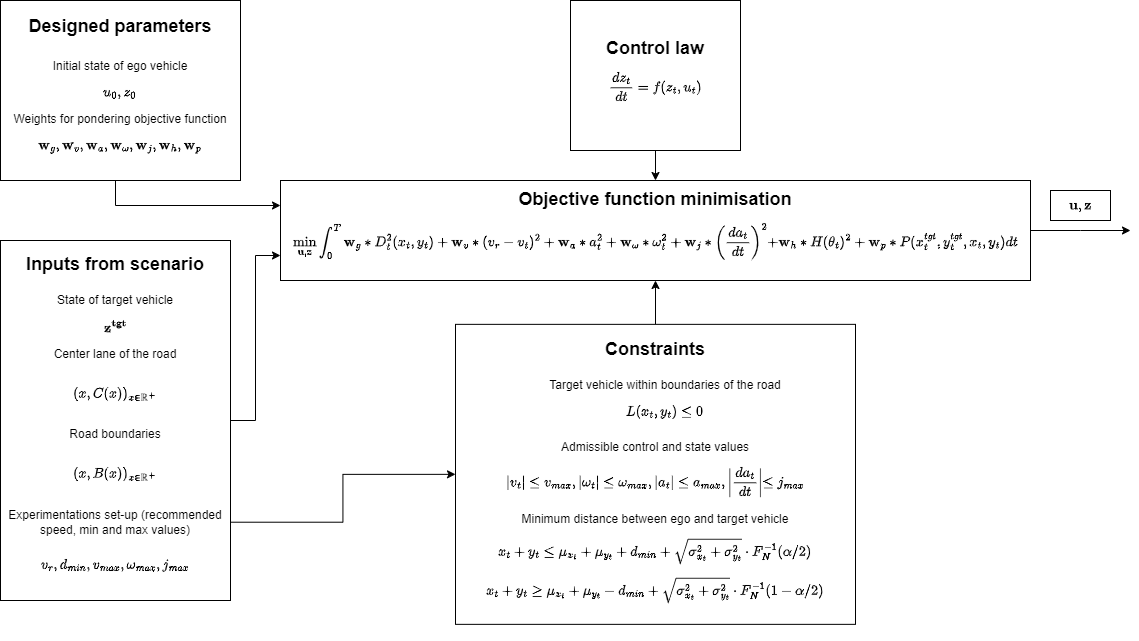}}
    \caption{Diagram of the optimal control under uncertainty.}
    \label{fig:control_diagram}
\end{figure}

Figure \ref{fig:control_diagram} is the control structure diagram of our optimal control problem. In trajectory planning, the goal is to find a control $\mathbf{u}$ over a time horizon $T$, so our ego's vehicle state $\mathbf{z}$ corresponds to a feasible solution for the autonomous vehicle.

\newpage

We formulated the problem of reference trajectory generation within the framework of autonomous vehicles and scenarios involving an ego and target vehicle. Section \ref{section:numerical_experiments} assesses the advantage of the continuous-time approach compared to discrete-time for various types of scenarios.

%% file: Numerical_Experiments.tex
\section{Numerical Experiments}
\label{section:numerical_experiments}

The numerical experiments aim to show the efficiency of the continuous-time model compared to previous discrete-time results \cite{zhang2022optimal}. The optimal control problem (\ref{eq:obj_1}) is solved by Python Package GEKKO \cite{beal2018gekko} in both deterministic and stochastic forms. We perform Model Predictive Control (MPC) \cite{rawlings2017model} which optimises the objective function and implicitly computes the model and its constraints at the same time.

The driving scenarios are sampled from the same experimental setup in discrete time to keep the comparison relevant. All scenarios explored in our numerical experiments in Subsection \ref{subsection:risky_scenario} and Subsection \ref{subsection:multiple_scenarios} were produced using SCANeR Studio software \cite{that2011integrated}. It is a tool designed for simulating the technology of Advanced Driver Assistance Systems (ADAS), made for assisting drivers.

All experimental results from this section are obtained by solving the trajectory planning problem for input scenarios from this software. Therefore, we assess the quality of our controller with respect to the industry's standards, with a high level of requirements to validate our method.

\begin{figure}[htp]
    \color{red}
    \centering
    \scalebox{.385}{\includegraphics{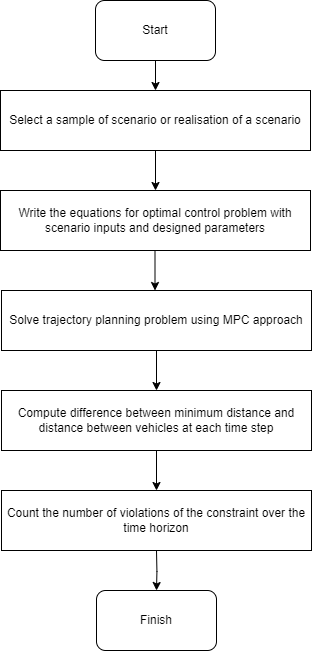}}
    \caption{Flowchart of the experiments' setup.}
    \label{fig:flowchart_experiments}
\end{figure}

Figure \ref{fig:flowchart_experiments} is a flowchart of the experiments and provides details about the metrics. As safety is our main concern, we assess the performance of our model by looking at constraint violations of the solution obtained. In both experiments from Subsection \ref{subsection:risky_scenario} and Subsection \ref{subsection:multiple_scenarios}, the process is the same. In the first experiment, we consider multiple realisations of the same scenario, that is, multiple measurements of the target's vehicle position. In the second experiment, multiple scenarios are considered.

We recall that this experiment setup is performed in an urban driving scenario with the following input parameters :

\begin{table}[ht]
\centering
\begin{tabular}{|c|c|c|}
\hline
Parameter                      & Function & Value \\ \hline
$v_r$                     &  Reference linear speed     &   $12\ m.s^{-1}$ \\ \hline
$d_{min}$   &  Minimum distance between ego and target vehicle  &   $5\ m$
\\ \hline
$v_{max}$   &    Maximum linear speed    &   $40\ m.s^{-1}$    \\ \hline
$\omega_{max}$  &         Maximum angular speed    &  $\frac{\pi}{6} \ s^{-1}$   \\ \hline
$j_{max}$       &     Maximum jerk    &  $0.6 \ m.s^{-3}$    \\ \hline
\end{tabular}
\caption{Parameters' values for urban driving scenarios during the simulation.}
\label{tab:parameters}
\end{table}

\begin{figure}[htp]
    \centering
    \scalebox{.5}{\includegraphics{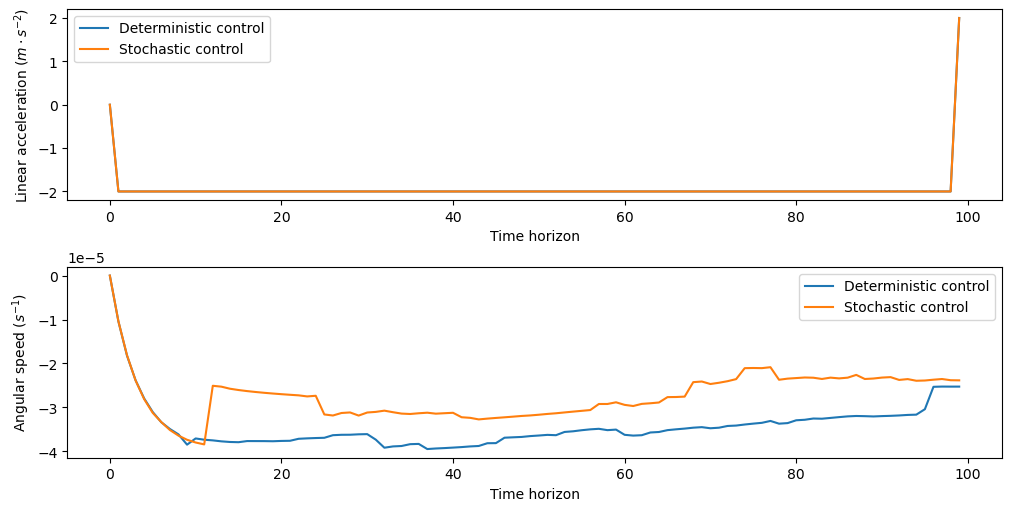}}
    \caption{Control values for an example scenario.}
    \label{fig:example}
\end{figure}

Figure \ref{fig:example} represents the control values (linear acceleration $(a_t)_{t \in [0,T]}$ and angular speed $(\omega_t)_{t \in [0,T]}$) in deterministic and stochastic controls for a given example scenario. We consider this scenario as risky, so the ego and target vehicles are willing to collide. This risky scenario is used for our first experiment in Section \ref{subsection:risky_scenario}.

After introducing the context in which we perform our experiments, we need to tackle the issues concerning the deterministic model. Subsection \ref{subsection:shortcomings_deterministic_model} depicts the limited results given by a purely deterministic model.

\subsection{Shortcomings of the deterministic model}
\label{subsection:shortcomings_deterministic_model}

\begin{figure}[htp]
    \centering
    \scalebox{.4}{\includegraphics{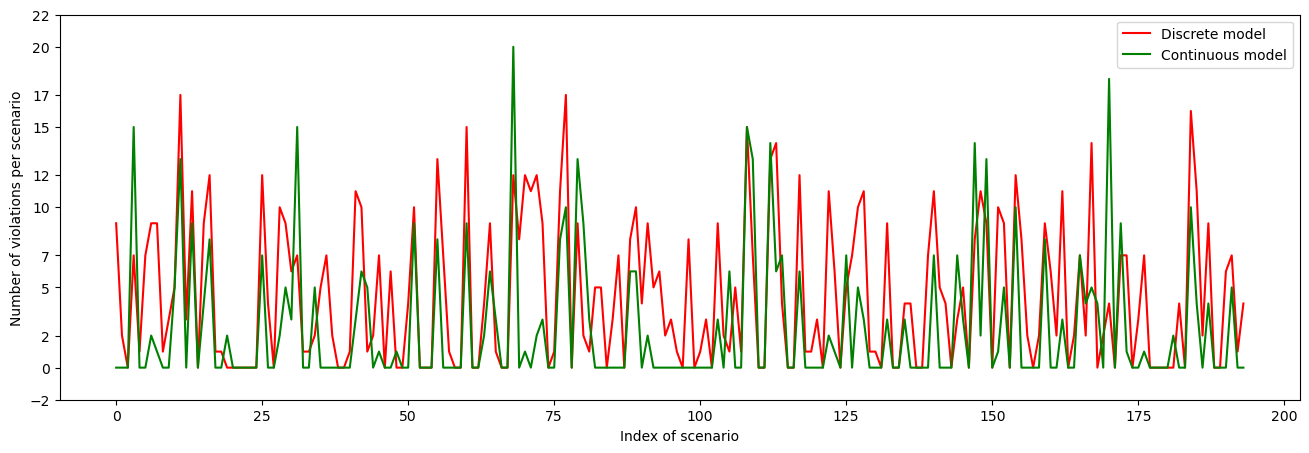}}
    \caption{Violations of constraints for deterministic models in discrete time and continuous time.}
    \label{fig:deterministic_violations}
\end{figure}

The deterministic model presents limited performances compared to the stochastic model and is subject to a high risk of collision in discrete time \cite{zhang2022optimal}. To show the robustness of the stochastic model, Figure \ref{fig:deterministic_violations} presents the number of violations for solutions given by the deterministic models in continuous time and discrete time on a corpus of 200 driving scenarios with different levels of risk.

\newpage

The number of violations of the constraint is very high for both models. We can conclude on eliminating the deterministic model in comparison with stochastic models. In the remainder of this study, we conduct our experiments for stochastic models only.

As we discard the use of deterministic models, we now focus on the design parameters of the model before performing the experiments. They are the weights $(\mathbf{w}_g, \mathbf{w}_v, \mathbf{w}_a, \mathbf{w}_{\omega}, \mathbf{w}_j, \mathbf{w}_h, \mathbf{w}_p)$ of the components of the objective function. Subsection \ref{subsection:design_weights} studies the impact of those weights to find a suitable solution.

\subsection{Impact of the continuous-time model and design parameters}
\label{subsection:design_weights}

\newpage

\begin{table*}[t]
\centering
\scriptsize
\begin{tabular}{clccccll}
\hline
\multicolumn{2}{|c|}{  $\mathbf{w}_g:\mathbf{w}_v:\mathbf{w}_a:\mathbf{w}_{\omega}:\mathbf{w}_j:\mathbf{w}_h:\mathbf{w}_p$} & \multicolumn{3}{c|}{5:1:1:1:1:1:1} & \multicolumn{3}{c|}{1:5:1:1:1:1:1}  \\ \hline
\multicolumn{2}{|c|}{$T$} & \multicolumn{1}{c|}{50} & \multicolumn{1}{c|}{200} & \multicolumn{1}{c|}{400} & \multicolumn{1}{c|}{50} & \multicolumn{1}{c|}{200} & \multicolumn{1}{c|}{400} \\ \hline
\multicolumn{2}{|c|}{Average CPU times} & \multicolumn{1}{c|}{0.261} & \multicolumn{1}{c|}{0.800} & \multicolumn{1}{c|}{1.079} & \multicolumn{1}{c|}{0.322} & \multicolumn{1}{c|}{0.794} & \multicolumn{1}{c|}{1.035} \\ \hline
\multicolumn{2}{|c|}{Average acceleration} & \multicolumn{1}{c|}{1.961} & \multicolumn{1}{c|}{1.98} & \multicolumn{1}{c|}{1.532} & \multicolumn{1}{c|}{1.849} & \multicolumn{1}{c|}{1.732} & \multicolumn{1}{c|}{1.373} \\ \hline
\multicolumn{2}{|c|}{Average angular velocity} & \multicolumn{1}{c|}{0} & \multicolumn{1}{c|}{0} & \multicolumn{1}{c|}{0} & \multicolumn{1}{c|}{0} & \multicolumn{1}{c|}{0} & \multicolumn{1}{c|}{0} \\ \hline
\multicolumn{2}{|c|}{Average velocity} & \multicolumn{1}{c|}{14.743} & \multicolumn{1}{c|}{15.741} & \multicolumn{1}{c|}{17.031} & \multicolumn{1}{c|}{14.718} & \multicolumn{1}{c|}{15.667} & \multicolumn{1}{c|}{16.861}  \\ \hline
\multicolumn{2}{|c|}{Average distance} & \multicolumn{1}{c|}{14.764} & \multicolumn{1}{c|}{31.523} & \multicolumn{1}{c|}{68.265} & \multicolumn{1}{c|}{14.739} & \multicolumn{1}{c|}{31.375} & \multicolumn{1}{c|}{67.589} \\ \hline
\multicolumn{1}{l}{} &  & \multicolumn{1}{l}{} & \multicolumn{1}{l}{} & \multicolumn{1}{l}{} & \multicolumn{1}{l}{} & 
 &  \\ \hline
\multicolumn{2}{|c|}{  $\mathbf{w}_g:\mathbf{w}_v:\mathbf{w}_a:\mathbf{w}_{\omega}:\mathbf{w}_j:\mathbf{w}_h:\mathbf{w}_p$} & \multicolumn{3}{c|}{1:1:5:1:1:1:1}  & \multicolumn{3}{c|}{1:1:1:5:1:1:1}  \\ \hline
\multicolumn{2}{|c|}{$T$} & \multicolumn{1}{c|}{50} & \multicolumn{1}{c|}{200} & \multicolumn{1}{c|}{400} & \multicolumn{1}{c|}{50} & \multicolumn{1}{c|}{200} & \multicolumn{1}{c|}{400} \\ \hline
\multicolumn{2}{|c|}{Average CPU times}  & \multicolumn{1}{c|}{0.247}   & \multicolumn{1}{c|}{0.744}   & \multicolumn{1}{c|}{1.145}   & \multicolumn{1}{c|}{0.251} & \multicolumn{1}{c|}{0.804}   & \multicolumn{1}{c|}{1.045}   \\ \hline
\multicolumn{2}{|c|}{Average acceleration} & \multicolumn{1}{c|}{1.960}   & \multicolumn{1}{c|}{1.972}   & \multicolumn{1}{c|}{1.437}   & \multicolumn{1}{c|}{1.960} & \multicolumn{1}{c|}{1.928}   & \multicolumn{1}{c|}{1.436}   \\ \hline
\multicolumn{2}{|c|}{Average angular velocity} & \multicolumn{1}{c|}{0}   & \multicolumn{1}{c|}{0}   & \multicolumn{1}{c|}{0}   & \multicolumn{1}{c|}{0}   & \multicolumn{1}{c|}{0}   & \multicolumn{1}{c|}{0}   \\ \hline
\multicolumn{2}{|c|}{Average velocity} & \multicolumn{1}{c|}{14.743}   & \multicolumn{1}{c|}{15.730}   & \multicolumn{1}{c|}{17.029}   & \multicolumn{1}{c|}{14.743}   & \multicolumn{1}{c|}{15.727}   & \multicolumn{1}{c|}{17.019}   \\ \hline
\multicolumn{2}{|c|}{Average distance} & \multicolumn{1}{c|}{14.764}   & \multicolumn{1}{c|}{31.501}   & \multicolumn{1}{c|}{68.256}   & \multicolumn{1}{c|}{14.764}   & \multicolumn{1}{c|}{31.496}   & \multicolumn{1}{c|}{68.220}   \\ \hline
\multicolumn{1}{l}{}               &               & \multicolumn{1}{l}{}    & \multicolumn{1}{l}{}    & \multicolumn{1}{l}{}    & \multicolumn{1}{l}{}    &                         &                         \\ \hline
\multicolumn{2}{|c|}{  $\mathbf{w}_g:\mathbf{w}_v:\mathbf{w}_a:\mathbf{w}_{\omega}:\mathbf{w}_j:\mathbf{w}_h:\mathbf{w}_p$} & \multicolumn{3}{c|}{1:1:1:1:5:1:1} & \multicolumn{3}{c|}{1:1:1:1:1:5:1}  \\ \hline
\multicolumn{2}{|c|}{$T$} & \multicolumn{1}{c|}{50} & \multicolumn{1}{c|}{200} & \multicolumn{1}{c|}{400} & \multicolumn{1}{c|}{50} & \multicolumn{1}{c|}{200} & \multicolumn{1}{c|}{400} \\ \hline
\multicolumn{2}{|c|}{Average CPU times} & \multicolumn{1}{c|}{0.245} & \multicolumn{1}{c|}{0.795} & \multicolumn{1}{c|}{1.044} & \multicolumn{1}{c|}{0.249} & \multicolumn{1}{c|}{0.947} & \multicolumn{1}{c|}{1.060} \\ \hline
\multicolumn{2}{|c|}{Average acceleration} & \multicolumn{1}{c|}{1.961} & \multicolumn{1}{c|}{1.928} & \multicolumn{1}{c|}{1.436} & \multicolumn{1}{c|}{1.960} & \multicolumn{1}{c|}{1.928} & \multicolumn{1}{c|}{1.436} \\ \hline
\multicolumn{2}{|c|}{Average angular velocity} & \multicolumn{1}{c|}{0} & \multicolumn{1}{c|}{0} & \multicolumn{1}{c|}{0} & \multicolumn{1}{c|}{0} & \multicolumn{1}{c|}{0} & \multicolumn{1}{c|}{0} \\ \hline
\multicolumn{2}{|c|}{Average velocity} & \multicolumn{1}{c|}{14.743} & \multicolumn{1}{c|}{15.727} & \multicolumn{1}{c|}{17.019} & \multicolumn{1}{c|}{14.743} & \multicolumn{1}{c|}{15.727} & \multicolumn{1}{c|}{17.019}  \\ \hline
\multicolumn{2}{|c|}{Average distance} & \multicolumn{1}{c|}{14.764} & \multicolumn{1}{c|}{31.496} & \multicolumn{1}{c|}{68.220} & \multicolumn{1}{c|}{14.739} & \multicolumn{1}{c|}{31.375} & \multicolumn{1}{c|}{68.220} \\ \hline
\multicolumn{1}{l}{} &  & \multicolumn{1}{l}{} & \multicolumn{1}{l}{} & \multicolumn{1}{l}{} & \multicolumn{1}{l}{} & 
 &  \\ \hline
\multicolumn{2}{|c|}{  $\mathbf{w}_g:\mathbf{w}_v:\mathbf{w}_a:\mathbf{w}_{\omega}:\mathbf{w}_j:\mathbf{w}_h:\mathbf{w}_p$} & \multicolumn{3}{c|}{1:1:1:1:1:1:5}  & \multicolumn{3}{c|}{1:1:1:1:1:1:1}  \\ \hline
\multicolumn{2}{|c|}{$T$} & \multicolumn{1}{c|}{50} & \multicolumn{1}{c|}{200} & \multicolumn{1}{c|}{400} & \multicolumn{1}{c|}{50} & \multicolumn{1}{c|}{200} & \multicolumn{1}{c|}{400} \\ \hline
\multicolumn{2}{|c|}{Average CPU times}  & \multicolumn{1}{c|}{0.242}   & \multicolumn{1}{c|}{0.799}   & \multicolumn{1}{c|}{1.038}   & \multicolumn{1}{c|}{0.241} & \multicolumn{1}{c|}{0.797}   & \multicolumn{1}{c|}{1.043}   \\ \hline
\multicolumn{2}{|c|}{Average acceleration} & \multicolumn{1}{c|}{1.960}   & \multicolumn{1}{c|}{1.928}   & \multicolumn{1}{c|}{1.436}   & \multicolumn{1}{c|}{1.960} & \multicolumn{1}{c|}{1.928}   & \multicolumn{1}{c|}{1.436}   \\ \hline
\multicolumn{2}{|c|}{Average angular velocity} & \multicolumn{1}{c|}{0}   & \multicolumn{1}{c|}{0}   & \multicolumn{1}{c|}{0}   & \multicolumn{1}{c|}{0}   & \multicolumn{1}{c|}{0}   & \multicolumn{1}{c|}{0}   \\ \hline
\multicolumn{2}{|c|}{Average velocity} & \multicolumn{1}{c|}{14.743} & \multicolumn{1}{c|}{15.727} & \multicolumn{1}{c|}{17.019} & \multicolumn{1}{c|}{14.743} & \multicolumn{1}{c|}{15.727} & \multicolumn{1}{c|}{17.019}  \\ \hline
\multicolumn{2}{|c|}{Average distance} & \multicolumn{1}{c|}{14.764} & \multicolumn{1}{c|}{31.496} & \multicolumn{1}{c|}{68.220} & \multicolumn{1}{c|}{14.739} & \multicolumn{1}{c|}{31.375} & \multicolumn{1}{c|}{68.220} \\ \hline
\multicolumn{1}{l}{}               &               & \multicolumn{1}{l}{}    & \multicolumn{1}{l}{}    & \multicolumn{1}{l}{}    & \multicolumn{1}{l}{}    &                         &                         \\ \hline
& \\ \hline
\end{tabular}
\caption{Comparison under different configurations over 100 different scenarios of the continuous-time model.}
\label{tab:compara}
\end{table*}

Table \ref{tab:compara} shows analytics for different configurations. As in \cite{zhang2022optimal}, the sensors are imperfect, and their measures contain errors which are modelled by adding a random noise with a normal distribution $\mathcal{N}(0,1)$ over the real positions of the target vehicle. To compensate the noise effect, we look at the average of control-state values over a sample of 100 different scenarios for each set of weights.

The time horizons $T$ considered here have higher values than the number of samples $N$ for discrete-time analysis \cite{zhang2022optimal}, thanks to dynamic solvers in continuous time, which can handle long term horizons. Steady-state real-time optimisation used for discrete-time optimal control problems requires solving $N$ different control-state equations, as each variable describes a time-step $k$. In contrast, dynamic solvers solve the equation for each control-state variable through the time horizon. In addition, the average CPU times are lower than for the discrete-time model \cite{zhang2022optimal}. The dynamics are similar for each set of weights, as they increase proportionally with the time horizons. The most significant increases in computational time are for weights $\mathbf{w}_a$ and $\mathbf{w}_v$. However, the continuous-time model is still faster and goes beyond the limits of the discrete-time model, as we develop in Section \ref{subsection:high_speeds}.

\begin{remark}
\label{remark:weights}
    The weights influence the behaviour of the ego vehicle to react more intensely to the weighted terms of the objective function. The most representative setup of a real-life simulation is with equal weights. We conducted the following experiments in Subsection \ref{subsection:risky_scenario} and Subsection \ref{subsection:multiple_scenarios} with this setup as imbalance introduces unwanted bias for the analysis of the results.
\end{remark}

\begin{figure}[htp]
    \centering
    \scalebox{.5}{\includegraphics{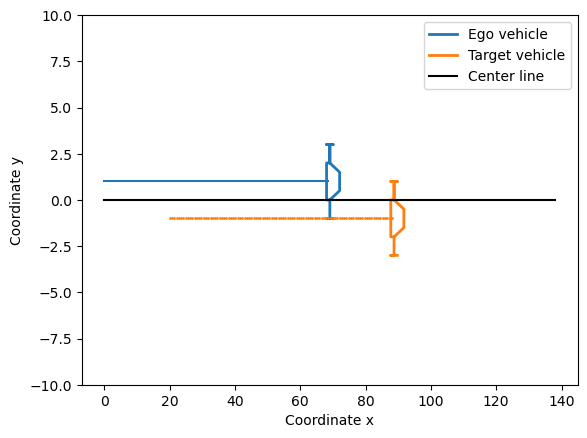}}
    \caption{Trajectory of ego and target vehicles.}
    \label{fig:trajectory}
\end{figure}

We now have all the theoretical elements to perform the simulations. The remainder of the study is focused on experiments. The first experiment from Subsection \ref{subsection:risky_scenario} illustrates the trajectory of a simple but risky scenario in Figure \ref{fig:trajectory}.

\subsection{Experiment 1: comparison of multiple realisations of one risky scenario}
\label{subsection:risky_scenario}

In the risky scenario, the recommended linear speed $v_r$ is $14\ m.s^{-1}$, which is $50.4\ km.h^{-1}$, while the maximum speed limit in cities is $50\ km.h^{-1}$.

Here, we suppose $(x^{tgt}_t, y^{tgt}_t)_{t \in [0, T]}$ to be sampled from random variables as in (\ref{eq:x_tgt_random}) and (\ref{eq:y_tgt_random}).

The number of steps for the discrete-time model and the time horizon for the continuous-time model are fixed with lower values for the first experiment than in the second experiment, where the scenarios are considered well-dimensioned. If they are not dimensioned correctly, the optimal control problem for trajectory planning does not guarantee the convergence to a feasible or realistic solution. As an example, if the linear recommended speed $v_r$ is too high and the associated weight $\mathbf{w}_v$ is large enough, the solution of the optimal control problem could give the ego vehicle to remain motionless at the beginning and wait for some time before accelerating. The horizon is lower to prevent this phenomenon in the simulation, especially in risky scenarios. When the speed increases, other parameters must be adapted to dimension different scenarios, as we do in Section \ref{subsection:high_speeds} to model road and highway driving scenarios.

\begin{figure}[htp]
    \centering
    \scalebox{.2}{\includegraphics{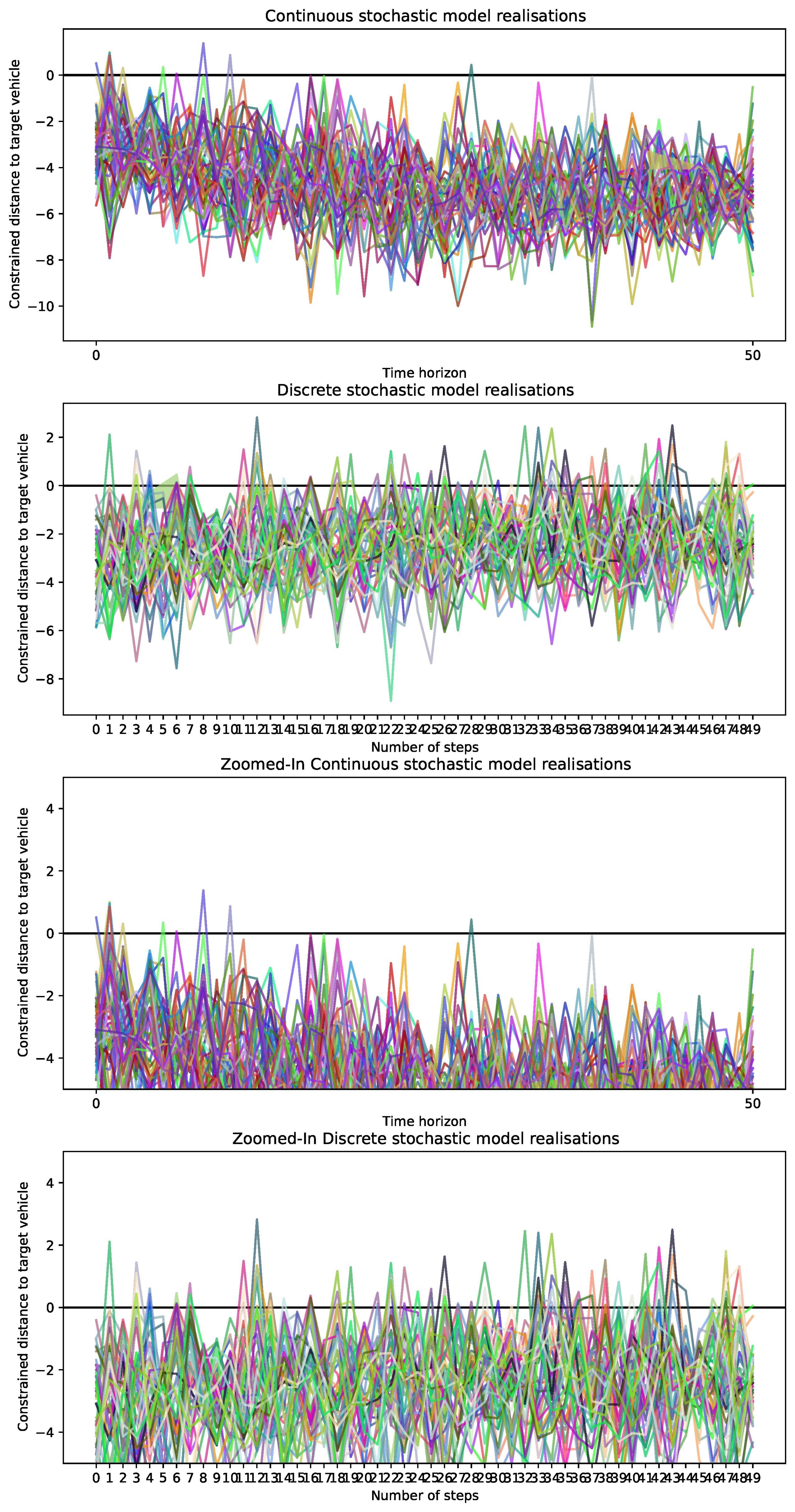}}
    \caption{Constraint function values of all realisations of the scenario for stochastic models.}
    \label{fig:one_scenario_plots}
\end{figure}

Figure \ref{fig:one_scenario_plots} shows the value of the quantity $d_{min} - K(x^{tgt}_t, y^{tgt}_t, x_t, y_t)$ derived from constraint (\ref{eq:ct7}) for each realisation, on the time interval $[0; T]$ for the continuous model and on the steps $\{1, ..., N\}$ for the discrete model. A constraint violation is visualised on positive values, which implies $d_{min} >  K(x^{tgt}_t, y^{tgt}_t, x_t, y_t)$, meaning vehicles are too close.

This chart presents results for both continuous-time and discrete-time stochastic solvers. The two plots on the bottom of the figure are zoomed-in figures of the two plots on the top. Zoomed-in plots are represented to visualise better the number of times the curve is above 0 during the simulation.\\
Each coloured curve is associated with a solution obtained for a realisation of the scenario. We performed 200 realisations, so 200 curves are represented. We simulated the same scenario 200 times and took 200 measurements of the position of the target's vehicle. During the time the curve is above the threshold of $0$, it means the vehicles are closer to each other than the minimum distance $d_{min}$ fixed in the constraint. This corresponds to a violation of the constraint.\\
On the y-axis, we visualise the distance between vehicles with respect to $d_{min}$. In continuous time, the ego vehicle stays further from the target vehicle than in discrete time. On the x-axis, we visualise the time horizon of the scenario. For the continuous-time model, the first moments of the simulation present more violations. Then, the model handles the ego vehicle to stay far enough from the target vehicle to respect the minimum distance. For the discrete-time model, violations are more frequent in the last moments, showing the model fails to control the vehicle not to violate the constraint.
\label{end_subsection:risky_scenario}

\newpage

\begin{figure}[htp]
    \centering
    \scalebox{.4}{\includegraphics{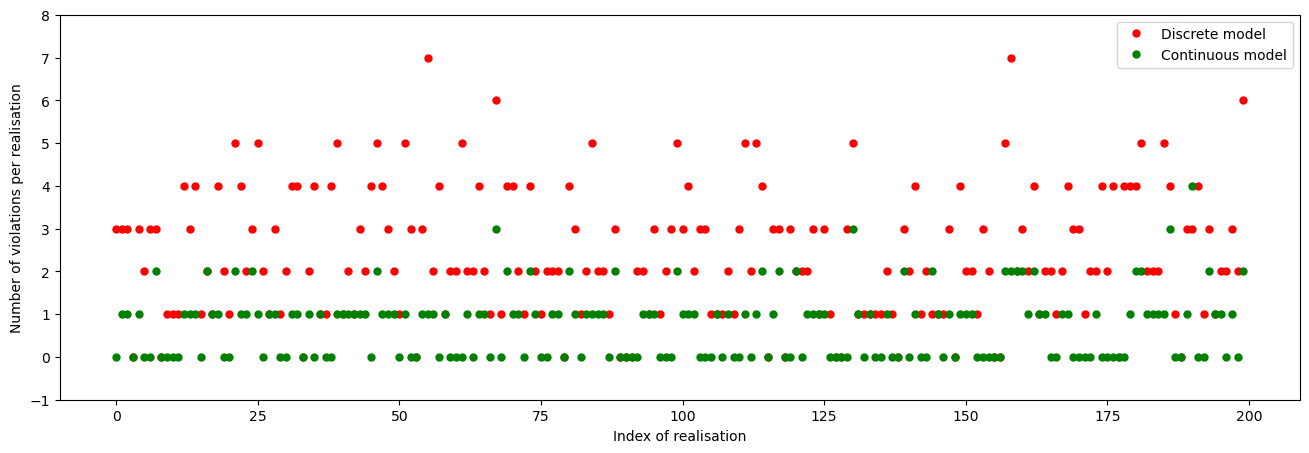}}
    \caption{Number of violations per realisation of the scenario for stochastic models.}
    \label{fig:one_scenario_points}
\end{figure}

Figure \ref{fig:one_scenario_points} represents the number of violations for each solution given by our model. Each solution is labelled for a realisation of the scenario, so there are 200 realisations represented on the x-axis. The number of violations is represented on the y-axis. Green dots are the number of violations in the solution given by the continuous-time model, and red dots are for solutions given by the discrete-time model. As we chose the scenario considered here as risky, we have more constraint violations in the realisations than for unstressed urban driving scenarios, as discussed in the second experiment Section \ref{subsection:multiple_scenarios}. However, the discrete-time model reaches the highest number of violations, and the average number of violations per realisation is higher for the discrete-time stochastic model than the continuous-time stochastic model.

\begin{figure}[htp]
    \centering
    \scalebox{.4}{\includegraphics{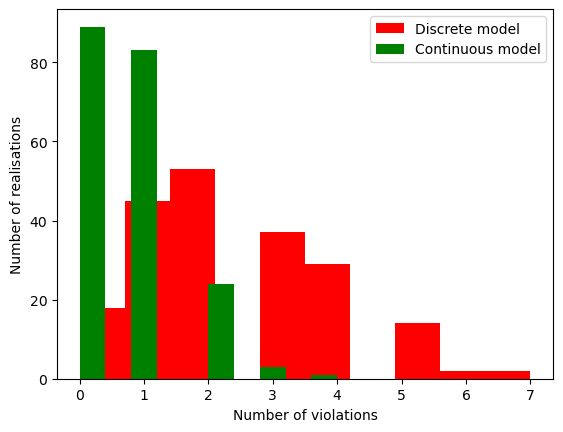}}
    \caption{Histogram of the number of violations per realisation of the scenario for stochastic models.}
    \label{fig:one_scenario_histogram}
\end{figure}

Figure \ref{fig:one_scenario_histogram} presents the histogram of the number of realisations per number of violations. We visualise that more than 160 solutions have one violation or less for continuous time. In contrast, more than half of solutions given by the discrete-time model have more than two violations.

\newpage

This first experiment presents the advantages of a continuous-time model compared to a discrete-time model in the specific case of a risky scenario. The following experiment in Subsection \ref{subsection:multiple_scenarios} extends to various classical urban driving scenarios.

\subsection{Experiment 2: comparison of multiple urban driving scenarios}
\label{subsection:multiple_scenarios}

\begin{figure}[p]
    \centering
    \scalebox{.5}{\includegraphics{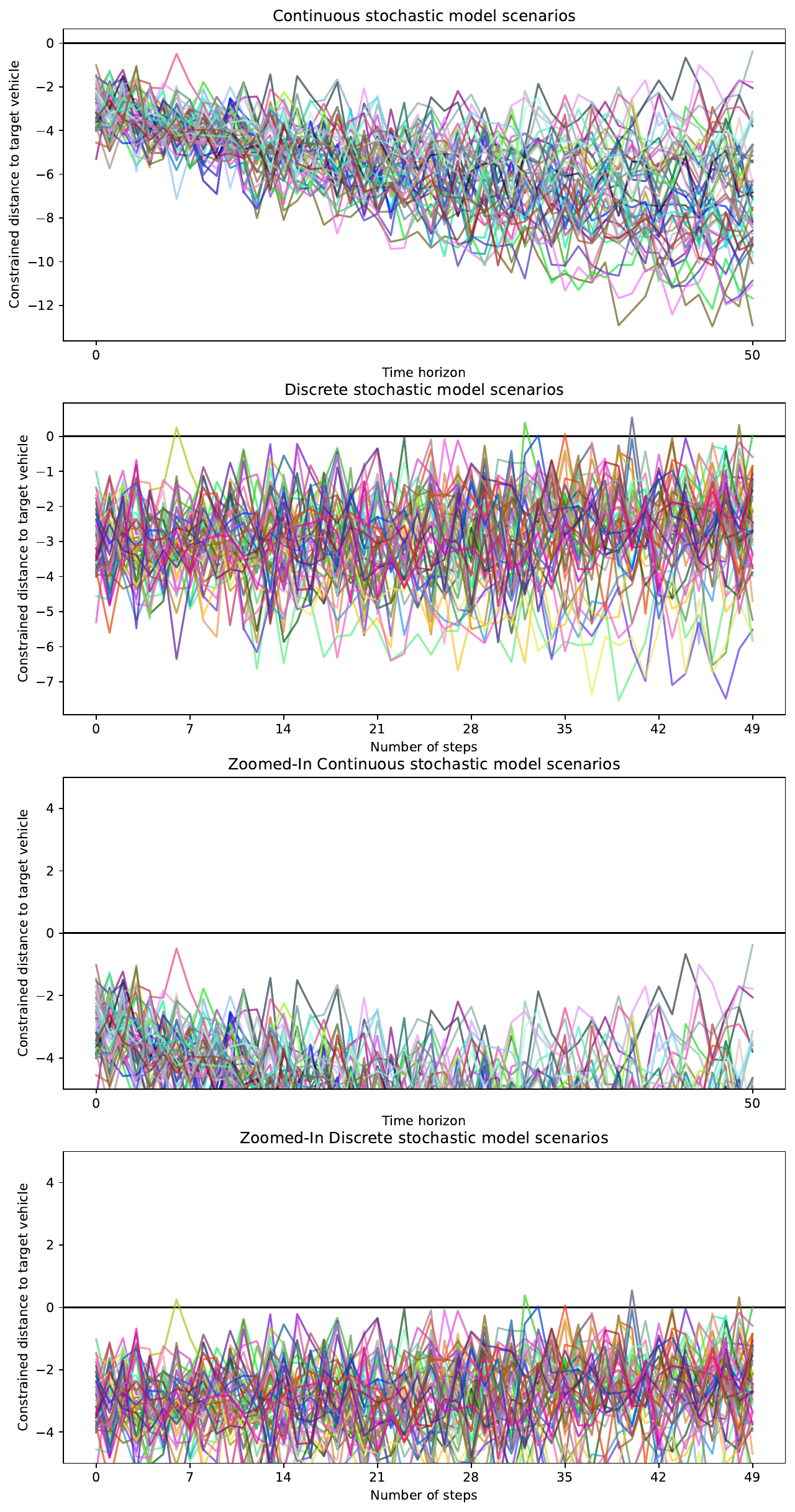}}
    \caption{Constraint function values of all scenarios for stochastic models.}
    \label{fig:multiple_scenarios_plots}
\end{figure}

\newpage

Figure \ref{fig:multiple_scenarios_plots} represent the same charts as the previous experiment in Subsection \ref{subsection:risky_scenario}, but each solution corresponds to a different scenario. Those scenarios are various urban driving scenarios generated from the simulator. As those are less risky scenarios than in the first experiment, both models perform better.

Each coloured curve is associated with a solution computed for a scenario. Solutions obtained in continuous time are more sparse as scenarios are more diverse compared to previous experiments, but results are better as the distance between vehicles is always respected. Over time, the risk of the vehicles colliding diminishes.\\
In discrete time, results are better than in the previous experiment, but they still present more violations than in continuous time. A higher risk of collision is observed as the number of steps of the simulation increases, meaning the model cannot handle the distance between vehicles above threshold $d_{min}$.

\begin{figure}[htp]
    \centering
    \scalebox{.4}{\includegraphics{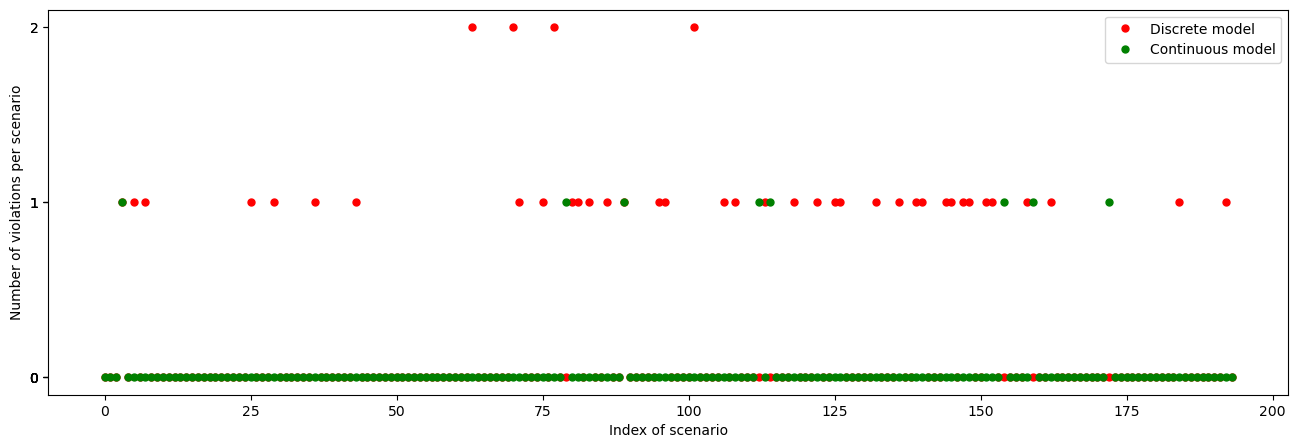}}
    \caption{Number of violations per scenario for stochastic models.}
    \label{fig:multiple_scenarios_points}
\end{figure}

\begin{figure}[htp]
    \centering
    \scalebox{.5}{\includegraphics{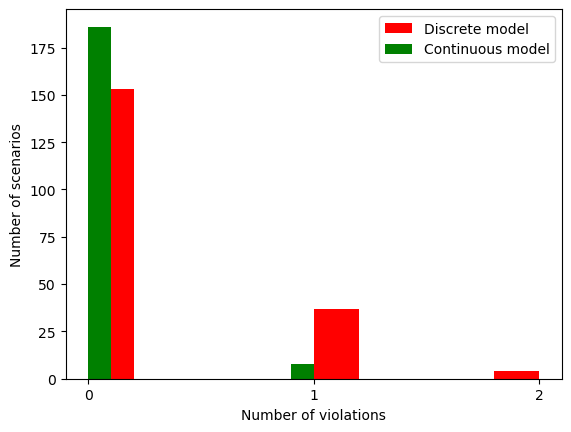}}
    \caption{Histogram of the number of violations per realisation of the scenario for stochastic models.}
    \label{fig:multiple_scenarios_histogram}
\end{figure}

Figure \ref{fig:multiple_scenarios_points} and \ref{fig:multiple_scenarios_histogram} are charts equivalent to Figure \ref{fig:one_scenario_points} and \ref{fig:one_scenario_histogram} in Subsection \ref{subsection:risky_scenario}. Observations show that in over 200 scenarios, the number of violations is small. The discrete-time model presents at most two violations and the continuous-time model presents at most one violation in a few scenarios. The histogram shows no violations for more than 180 scenarios for the continuous-time model and more than 150 scenarios for the discrete-time one.

\newpage

Our experiments proved the efficiency of our method. To assess experimentally the effectiveness of the constraints derived in (\ref{eq:ct8}) and (\ref{eq:ct9}), following Subsection \ref{subsection:chance_constraints} presents the ratio of time the minimum distance between vehicles is respected with respect to the threshold $\alpha$.

\newpage
        
\subsection{Chance-constraints in stochastic models}
\label{subsection:chance_constraints}
\begin{figure}[htp]
    \centering
    \scalebox{.5}{\includegraphics{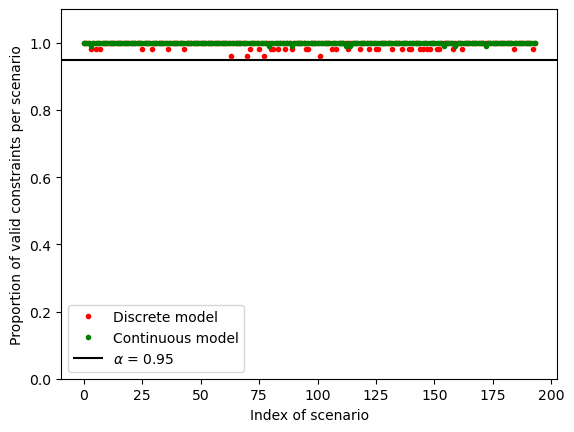}}
    \caption{Proportion of valid constraints per scenario for discrete-time and continuous-time stochastic models}
    \label{fig:multiple_scenario_alpha_threshold}
\end{figure}

Figure \ref{fig:multiple_scenario_alpha_threshold} is a representation of the chance constraint (\ref{eq:proba_alpha}) over all scenarios for both discrete-time and continuous-time models. We set a confidence value $\alpha = 0.95$, represented as a bold line on this figure. The deterministic equivalent second-order conic problem constraints (\ref{eq:x_tgt_random}) and (\ref{eq:y_tgt_random}) are sufficiently robust to guarantee the ego vehicle does not violate the constraint $95\%$ of the time through one simulation. The continuous-time stochastic model performs better than the discrete-time stochastic model, as for 200 scenarios, the model is valid more than $99\%$ of the time.

We showed that the continuous-time model performs better for urban driving scenarios, but the benefit compared to the discrete-time approach is not high, as both models perform well. Subsection \ref{subsection:high_speeds} extend solvers to high-speed scenarios on national roads and highways where the discrete-time model often fails to find a solution.

\subsection{Robustness of continuous-time model to long term horizons $T$}
\label{subsection:high_speeds}

The advantage of the continuous-time model is to be able to handle a long term horizon ($T \geq 100$) where the discrete model presents limitations over the number of steps $N$. Therefore, the continuous model can handle scenarios other than those studied in Subsection \ref{subsection:risky_scenario} and Subsection \ref{subsection:multiple_scenarios}, allowing higher speeds such as driving on roads outside city restrictions and highways. This approach requires more data points but is more robust to prevent unexpected events and have a precautionary approach to driving.

\begin{table}[ht]
\begin{center}
\begin{tabular}{|c|c|c|c|c|c|c|}
    \cline{2-7}
    \multicolumn{1}{c|}{} & \multicolumn{3}{c|}{Continuous model} & \multicolumn{3}{c|}{Discrete model} \\ \hline
    $v_r \slash \omega_{max}$ & $\frac{\pi}{6}$ & $\frac{\pi}{4}$ & $\frac{\pi}{2}$ & $\frac{\pi}{6}$ & $\frac{\pi}{4}$ & $\frac{\pi}{2}$\\ \hline
    $22\ m.s^{-1}$ & $96$ & $98$ & $99$ & $6$ & $7$ & $2$\\ \hline
    $36\ m.s^{-1}$ & $85$ & $95$ & $92$ & $3$ & $4$ & $1$\\ \hline
\end{tabular}
\end{center}
\caption{Number of solutions found over 100 scenarios}
\label{table:high_speeds}
\end{table}

Table \ref{table:high_speeds} presents the number of solutions found for 100 scenarios, given different constraints over the control. For this experiment, the computation time is, on average, $30$ minutes for the discrete-time model and $5$ minutes for the continuous-time one.

The discrete-time model does not find a solution for all real-life scenarios. The steady-state solver is limited in this case. The discretisation of the problem does not reflect well real-life driving and results in a loss of information. The dynamic solvers is robust enough to find feasible solutions in continuous time when the scenarios are more diverse concerning target vehicle speed and centre lane variations.

Our model showed promising results for a broad range of scenarios, meaning the solution obtained is feasible and allows the ego vehicle to perform safe driving for most simulations. We developed a generic model that can be extended to other cases.

\section{Limitations of the model and future research directions}
\label{section:limitations}

In this study, our experiments proved the interest of a continuous-time approach using dynamic solvers for solving the optimal control problem of trajectory planning for autonomous vehicles. Experiments were conducted to validate the stochastic approach for handling constraints on the target's vehicle position. We presented a generic model that can be derived for other vehicles. Therefore, one can develop a model for other environments. Our model does not consider how to overtake a car and the associated constraints of such a scenario. It corresponds to an environment other than the one studied here, so it must be the subject of future research to address this problem. Other types of vehicles also embed constraints from their environment. For spacecraft, aircraft or submarines, each environment requires additional constraints than those considered in our generic model to perform control of those vehicles. Physical parameters of height, depth, wind velocity or strength of ocean currents are mandatory to include in the control, which opens opportunities for further research.\\
From a theoretical perspective, different probability distributions are studied in the literature on chance constraints. Recent researches discuss about elliptical distributions \cite{nguyen2022random} and mixture distributions \cite{peng2021chance} to model the stochasticity of the chance constraints. The application of chance constraints to the optimal control problem of trajectory planning can lead to the exploration of other distributions when the environment is more complex.

%% file: Conclusion.tex
\section{Conclusions and practical recommendations}
\label{section:conclusion}

This paper shows the effectiveness of a continuous-time approach for optimal control-based reference trajectory generation compared to a discrete-time framework. Chance-constrained optimisation addresses the stochastic aspect of the model, which considers the uncertainty of information concerning autonomous vehicles. Our analysis showed low performance for deterministic models compared to the stochastic models in continuous time. Promising results are obtained for urban driving scenarios with both continuous-time and discrete-time stochastic models, with slightly better results in continuous time, specifically for risky scenarios where it diminishes the risks of collisions. Still, the discrete-time stochastic model fails to find solutions for national roads and highways with higher speed limits, as those scenarios require planning over long term horizons. The computation time is higher for the discrete-time stochastic model as the number of control-state variables increases with the horizon. The continuous-time stochastic model reflects real-life driving conditions better and converges faster to an optimal solution. The model developed shows its effectiveness in using the Model Predictive Control approach for solving a trajectory planning problem. The robustness of the chance constraints approach gives a model with a controlled level of risk.

Chance-constraints approach gives perspectives to control various types of autonomous vehicles, as discussed in Section \ref{section:limitations}. To model a complex environment accurately, further research should investigate both recent research in chance constraints and advances in control theory and its application, as those mentioned in introduction Section \ref{section:introduction} and Section \ref{section:related_work}.

%% file: main.bbl
\begin{thebibliography}{10}

\bibitem{zhang2022optimal}
Shangyuan Zhang, Makhlouf Hadji, and Abdel Lisser.
\newblock Optimal control based trajectory planning under uncertainty.
\newblock In {\em International Conference on Intelligent Transport Systems}, pages 73--88. Springer, 2022.

\bibitem{guo2023survey}
Yuqing Guo, Zelin Guo, Yazhou Wang, Danya Yao, Bai Li, and Li~Li.
\newblock A survey of trajectory planning methods for autonomous driving—part i: Unstructured scenarios.
\newblock {\em IEEE Transactions on Intelligent Vehicles}, 2023.

\bibitem{li2023real}
Guoqiang Li, Xudong Zhang, Hongliang Guo, Basilio Lenzo, and Ningyuan Guo.
\newblock Real-time optimal trajectory planning for autonomous driving with collision avoidance using convex optimization.
\newblock {\em Automotive Innovation}, 6(3):481--491, 2023.

\bibitem{takemura2024uncertainty}
Reiya Takemura and Genya Ishigami.
\newblock Uncertainty-aware trajectory planning: Using uncertainty quantification and propagation in traversability prediction of planetary rovers.
\newblock {\em IEEE Robotics \& Automation Magazine}, 2024.

\bibitem{le2024stochastic}
Viet-Anh Le, Behdad Chalaki, Filippos~N Tzortzoglou, and Andreas~A Malikopoulos.
\newblock Stochastic time-optimal trajectory planning for connected and automated vehicles in mixed-traffic merging scenarios.
\newblock {\em IEEE Transactions on Control Systems Technology}, 2024.

\bibitem{hentout2023review}
Abdelfetah Hentout, Abderraouf Maoudj, and Mustapha Aouache.
\newblock A review of the literature on fuzzy-logic approaches for collision-free path planning of manipulator robots.
\newblock {\em Artificial Intelligence Review}, 56(4):3369--3444, 2023.

\bibitem{chu2024motion}
Wenbo Chu, Kai Yang, Shen Li, and Xiaolin Tang.
\newblock Motion planning for autonomous driving with real traffic data validation.
\newblock {\em Chinese Journal of Mechanical Engineering}, 37(1):6, 2024.

\bibitem{anderson2010optimal}
Sterling~J Anderson, Steven~C Peters, Tom~E Pilutti, and Karl Iagnemma.
\newblock An optimal-control-based framework for trajectory planning, threat assessment, and semi-autonomous control of passenger vehicles in hazard avoidance scenarios.
\newblock {\em International Journal of Vehicle Autonomous Systems}, 8(2-4):190--216, 2010.

\bibitem{choi2023motion}
Deok-Kee Choi.
\newblock Motion tracking of four-wheeled mobile robots in outdoor environments using bayes’ filters.
\newblock {\em International Journal of Precision Engineering and Manufacturing}, 24(5):767--786, 2023.

\bibitem{sobanski2024predefined}
Rafa{\l}~M Soba{\'n}ski, Maciej~M Micha{\l}ek, and Michael Defoort.
\newblock Predefined-time vfo control design for unicycle-like mobile robots.
\newblock {\em Nonlinear Dynamics}, 112(5):3591--3603, 2024.

\bibitem{adzkiya2024control}
Dieky Adzkiya, Muhammad~Syifa'ul Mufid, Febrianti~Silviana Saputri, and Alessandro Abate.
\newblock Control design of discrete-time unicycle model using satisfiability modulo theory.
\newblock {\em Systems Science \& Control Engineering}, 12(1):2316166, 2024.

\bibitem{degorre2023survey}
Lo{\"\i}ck Degorre, Emmanuel Delaleau, and Olivier Chocron.
\newblock A survey on model-based control and guidance principles for autonomous marine vehicles.
\newblock {\em Journal of Marine Science and Engineering}, 11(2):430, 2023.

\bibitem{vishnu2023improving}
Chalavadi Vishnu, Vineel Abhinav, Debaditya Roy, C~Krishna Mohan, and Ch~Sobhan Babu.
\newblock Improving multi-agent trajectory prediction using traffic states on interactive driving scenarios.
\newblock {\em IEEE Robotics and Automation Letters}, 8(5):2708--2715, 2023.

\bibitem{zhang2022adversarial}
Qingzhao Zhang, Shengtuo Hu, Jiachen Sun, Qi~Alfred Chen, and Z~Morley Mao.
\newblock On adversarial robustness of trajectory prediction for autonomous vehicles.
\newblock In {\em Proceedings of the IEEE/CVF Conference on Computer Vision and Pattern Recognition}, pages 15159--15168, 2022.

\bibitem{ren2022chance}
Kai Ren, Heejin Ahn, and Maryam Kamgarpour.
\newblock Chance-constrained trajectory planning with multimodal environmental uncertainty.
\newblock {\em IEEE Control Systems Letters}, 7:13--18, 2022.

\bibitem{vaskov2024friction}
Sean Vaskov, Rien Quirynen, Marcel Menner, and Karl Berntorp.
\newblock Friction-adaptive stochastic nonlinear model predictive control for autonomous vehicles.
\newblock {\em Vehicle System Dynamics}, 62(2):347--371, 2024.

\bibitem{wang2020non}
Allen Wang, Ashkan Jasour, and Brian~C Williams.
\newblock Non-gaussian chance-constrained trajectory planning for autonomous vehicles under agent uncertainty.
\newblock {\em IEEE Robotics and Automation Letters}, 5(4):6041--6048, 2020.

\bibitem{schwarting2018planning}
Wilko Schwarting, Javier Alonso-Mora, and Daniela Rus.
\newblock Planning and decision-making for autonomous vehicles.
\newblock {\em Annual Review of Control, Robotics, and Autonomous Systems}, 1:187--210, 2018.

\bibitem{rigatos2020nonlinear}
G~Rigatos, P~Siano, F~Zouari, and Sul Ademi.
\newblock Nonlinear optimal control of autonomous submarines’ diving.
\newblock {\em Marine Systems \& Ocean Technology}, 15:57--69, 2020.

\bibitem{rigatos2024nonlinear}
G~Rigatos, J~Pomares, M~Abbaszadeh, K~Busawon, Z~Gao, F~Zouari, and Ecole~d'Ing enieurs~de Tunis.
\newblock Nonlinear optimal control for free-floating space robotic manipulators.
\newblock {\em Spacecraft Satellite}, 1:563, 2024.

\bibitem{zouari2024finite}
Farouk Zouari, Asier Ibeas, Abdesselem Boulkroune, and Jinde Cao.
\newblock Finite-time adaptive event-triggered output feedback intelligent control for noninteger order nonstrict feedback systems with asymmetric time-varying pseudo-state constraints and nonsmooth input nonlinearities.
\newblock {\em Communications in Nonlinear Science and Numerical Simulation}, 136:108036, 2024.

\bibitem{zouari2018neural}
Farouk Zouari and Amina Boubellouta.
\newblock Neural approximation-based adaptive control for pure-feedback fractional-order systems with output constraints and actuator nonlinearities.
\newblock In {\em Advanced Synchronization Control and Bifurcation of Chaotic Fractional-Order Systems}, pages 468--495. IGI Global, 2018.

\bibitem{wu2023ccgnet}
Dawen Wu and Abdel Lisser.
\newblock Ccgnet: A deep learning approach to predict nash equilibrium of chance-constrained games.
\newblock {\em Information Sciences}, 627:20--33, 2023.

\bibitem{chen2022deep}
Lienhung Chen, Zhongliang Jiang, Long Cheng, Alois~C Knoll, and Mingchuan Zhou.
\newblock Deep reinforcement learning based trajectory planning under uncertain constraints.
\newblock {\em Frontiers in Neurorobotics}, 16:883562, 2022.

\bibitem{peng2021chance}
Shen Peng, Navnit Yadav, Abdel Lisser, and Vikas~Vikram Singh.
\newblock Chance-constrained games with mixture distributions.
\newblock {\em Mathematical Methods of Operations Research}, 94:71--97, 2021.

\bibitem{tassouli2022solving}
Siham Tassouli and Abdel Lisser.
\newblock Solving linear programs with joint probabilistic constraints with dependent rows using a dynamical neural network.
\newblock {\em Results in Control and Optimization}, 9:100178, 2022.

\bibitem{liu2022distributionally}
Jia Liu, Abdel Lisser, and Zhiping Chen.
\newblock Distributionally robust chance constrained geometric optimization.
\newblock {\em Mathematics of Operations Research}, 47(4):2950--2988, 2022.

\bibitem{varagapriya2023joint}
V~Varagapriya, Vikas~Vikram Singh, and Abdel Lisser.
\newblock Joint chance-constrained markov decision processes.
\newblock {\em Annals of Operations Research}, 322(2):1013--1035, 2023.

\bibitem{peng2022bounds}
Shen Peng, Francesca Maggioni, and Abdel Lisser.
\newblock Bounds for probabilistic programming with application to a blend planning problem.
\newblock {\em European Journal of Operational Research}, 297(3):964--976, 2022.

\bibitem{ornelas2010discrete}
Fernando Ornelas, Edgar~N Sanchez, and Alexander~G Loukianov.
\newblock Discrete-time inverse optimal control for nonlinear systems trajectory tracking.
\newblock In {\em 49th IEEE Conference on Decision and Control (CDC)}, pages 4813--4818. IEEE, 2010.

\bibitem{liu2017path}
Chang Liu, Seungho Lee, Scott Varnhagen, and H~Eric Tseng.
\newblock Path planning for autonomous vehicles using model predictive control.
\newblock In {\em 2017 IEEE Intelligent Vehicles Symposium (IV)}, pages 174--179. IEEE, 2017.

\bibitem{prekopa2013stochastic}
Andr{\'a}s Pr{\'e}kopa.
\newblock {\em Stochastic programming}, volume 324.
\newblock Springer Science \& Business Media, 2013.

\bibitem{beal2018gekko}
Logan Beal, Daniel Hill, R~Martin, and John Hedengren.
\newblock Gekko optimization suite.
\newblock {\em Processes}, 6(8):106, 2018.

\bibitem{rawlings2017model}
James~Blake Rawlings, David~Q Mayne, and Moritz Diehl.
\newblock {\em Model predictive control: theory, computation, and design}, volume~2.
\newblock Nob Hill Publishing Madison, WI, 2017.

\bibitem{that2011integrated}
Thomas~Nguen That and Jordi Casas.
\newblock An integrated framework combining a traffic simulator and a driving simulator.
\newblock {\em Procedia-Social and Behavioral Sciences}, 20:648--655, 2011.

\bibitem{nguyen2022random}
Hoang~Nam Nguyen, Abdel Lisser, and Vikas~Vikram Singh.
\newblock Random games under elliptically distributed dependent joint chance constraints.
\newblock {\em Journal of Optimization Theory and Applications}, 195(1):249--264, 2022.

\end{thebibliography}
